\documentclass[10pt]{article}
\usepackage{cite}
\usepackage{mathrsfs}
\usepackage{amsfonts}
\usepackage{amsmath}
\usepackage{amsfonts,amssymb}
\usepackage{dsfont}
\usepackage{curves}
\usepackage{mathrsfs}
\usepackage{pifont}
\usepackage{amssymb}
\allowdisplaybreaks

\numberwithin{equation}{section}

\date{}

\def\BigRoman{\uppercase\expandafter{\romannumeral\number\count 255 }}
\def\Romannumeral{\afterassignment\BigRoman\count255=}

\begin{document}
\title{An $A_{\alpha}$-spectral radius for the existence of $\{P_3,P_4,P_5\}$-factors in graphs
}
\author{\small  Yuli Zhang$^{1}$, Sizhong Zhou$^{2}$\footnote{Corresponding
author. E-mail address: zsz\_cumt@163.com (S. Zhou)}\\
\small $1$. School of Science, Dalian Jiaotong University,\\
\small Dalian, Liaoning 116028, China\\
\small $2$. School of Science, Jiangsu University of Science and Technology,\\
\small Zhenjiang, Jiangsu 212100, China\\
}

\maketitle
\begin{abstract}
\noindent Let $G$ be a connected graph of order $n$ with $n\geq25$. A $\{P_3,P_4,P_5\}$-factor is a spanning subgraph $H$ of $G$ such that every
component of $H$ is isomorphic to an element of $\{P_3,P_4,P_5\}$. Nikiforov introduced the $A_{\alpha}$-matrix of $G$ as
$A_{\alpha}(G)=\alpha D(G)+(1-\alpha)A(G)$ [V. Nikiforov, Merging the $A$- and $Q$-spectral theories, Appl. Anal. Discrete Math. 11 (2017) 81--107],
where $\alpha\in[0,1]$, $D(G)$ denotes the diagonal matrix of vertex degrees of $G$ and $A(G)$ denotes the adjacency matrix of $G$. The largest
eigenvalue of $A_{\alpha}(G)$, denoted by $\lambda_{\alpha}(G)$, is called the $A_{\alpha}$-spectral radius of $G$. In this paper, it is proved
that $G$ has a $\{P_3,P_4,P_5\}$-factor unless $G=K_1\vee(K_{n-2}\cup K_1)$ if $\lambda_{\alpha}(G)\geq\lambda_{\alpha}(K_1\vee(K_{n-2}\cup K_1))$,
where $\alpha$ be a real number with $0\leq\alpha<\frac{2}{3}$.
\\
\begin{flushleft}

(2020) Mathematics Subject Classification: 05C50, 05C70, 05C38
\end{flushleft}
\end{abstract}

\section{Introduction}

We deal with finite undirected graphs without loops or multiple edges. Let $G$ denote a graph with vertex set $V(G)$ and edge set $E(G)$. For a vertex
$v\in V(G)$, the neighborhood of $v$ and the degree of $v$ in $G$ are denoted by $N_G(v)$ and $d_G(v)$, respectively. Let $i(G)$ denote the number of
isolated vertices in $G$. For a subset $S\subseteq V(G)$, let $G[S]$ and $G-S$ denote the subgraphs of $G$ induced by $S$ and $V(G)-S$, respectively.
For two vertex disjoint graphs $G_1$ and $G_2$, the union of $G_1$ and $G_2$ is denoted by $G_1\cup G_2$. Let $tG$ stand for the disjoint union of $t$
copies of $G$, where $t$ is a positive integer. The join $G_1\vee G_2$ is the graph obtained by joining each vertex of $G_1$ to each vertex of $G_2$.
We denote the path, the cycle, the star and the complete graph of order $n$ by $P_n$, $C_n$, $K_{1,n-1}$ and $K_n$, respectively. Let $c$ be a real
number. Recall that $\lfloor c\rfloor$ is the greatest integer with $\lfloor c\rfloor\leq c$.

Let $\mathcal{H}$ denote a set of connected graphs. Then a spanning subgraph $H$ of $G$ is called an $\mathcal{H}$-factor if every component of $H$ is
an element of $\mathcal{H}$. If $\mathcal{H}=\{P_3,P_4,P_5\}$, then an $\mathcal{H}$-factor is called a $\{P_3,P_4,P_5\}$-factor. Write
$P_{\geq k}=\{P_i|i\geq k\}$. If $\mathcal{H}=P_{\geq k}$, then an $\mathcal{H}$-factor is called a $P_{\geq k}$-factor. If $\mathcal{H}=\{K_2,C_i|i\geq3\}$,
then an $\mathcal{H}$-factor is called a $\{K_2,C_i|i\geq3\}$-factor. If $\mathcal{H}=\{K_{1,j}|1\leq j\leq k\}$, then an $\mathcal{H}$-factor is
called a $\{K_{1,j}|1\leq j\leq k\}$-factor.

Kano, Lu and Yu \cite{KLY} established a connection between the number of isolated vertices and $\{P_3,P_4,P_5\}$-factors in graphs. Akiyama, Avis and
Era \cite{AAE} proved that a graph $G$ contains a $P_{\geq2}$-factor if and only if $i(G-S)\leq2|S|$ for any subset $S\subseteq V(G)$. Kaneko \cite{Ka}
provided a characterization of a graph having a $P_{\geq3}$-factor. Liu and Pan \cite{LP2}, Gao, Wang and Chen \cite{GWC}, Dai and Hu \cite{DH}, Zhou
et al. \cite{ZSL,Zs,Zp,ZSB} obtained some sufficient conditions on the existence of $P_{\geq2}$-factors and $P_{\geq3}$-factors in graphs. Tutte \cite{Tt}
got a criterion for a graph containing a $\{K_2,C_i|i\geq3\}$-factor. Klopp and Steffen \cite{KS} investigated the properties of
$\{K_{1,1},K_{1,2},C_i|i\geq3\}$-factors in graphs. Amahashi and Kano \cite{AK} posed a criterion for a graph with a $\{K_{1,j}|1\leq j\leq k\}$-factor,
where $k$ is an integer with $k\geq2$. Zhou, Xu and Sun \cite{ZXS} studied the existence of a $\{K_{1,j}|1\leq j\leq k\}$-factor in a graph, where $k$
is an integer with $k\geq2$. Kano and Saito \cite{KSs} showed a sufficient condition for a graph to contain a $\{K_{1,j}|k\leq j\leq2k\}$-factor, where
$k$ is an integer with $k\geq2$. Zhou, Bian and Sun \cite{ZBS} obtained two sufficient conditions for a graph $G$ with $\delta(G)\geq2$ to have a
$\{K_{1,j}|1\leq j\leq k,\mathcal{T}(2k+1)\}$-factor, where $k\geq2$ is an integer and $\mathcal{T}(2k+1)$ is a special class of trees. For many other
results on spanning subgraphs, we refer the readers to \cite{GWW,Za1,ZPX1,ZPX}.

Let $A(G)$ and $D(G)$ denote the adjacency matrix and the degree diagonal matrix of $G$, respectively. We use $\lambda(G)$ to denote the adjacency
spectral radius of $G$. Let $Q(G)=D(G)+A(G)$ be the signless Laplacian matrix of $G$. The signless Laplacian spectral radius of $G$ is denoted by $q(G)$.
For any $\alpha\in[0,1)$, Nikiforov \cite{N} introduced the $A_{\alpha}$-matrix of $G$ as
$$
A_{\alpha}(G)=\alpha D(G)+(1-\alpha)A(G).
$$
Notice that $A_{\alpha}(G)=A(G)$ if $\alpha=0$ and $A_{\alpha}(G)=\frac{1}{2}Q(G)$ if $\alpha=\frac{1}{2}$. The largest eigenvalue of $A_{\alpha}(G)$,
denoted by $\lambda_{\alpha}(G)$, is called the $A_{\alpha}$-spectral radius of $G$. Clearly, $\lambda_0(G)$ is the adjacency spectral radius of $G$ and
$2\lambda_\frac{1}{2}(G)$ is the signless Laplacian spectral radius of $G$. Thus, $\lambda_{\alpha}(G)$ generalizes both the adjacency spectral radius
and the signless Laplacian spectral radius of $G$. In recent years, the $A_{\alpha}$-matrix of $G$ attracts a great deal of attention. For details,
we refer the readers to \cite{NR,ABGD,LYL,Wc,ZSL1,ZZL}.

O \cite{O}, Zhao, Huang and Wang \cite{ZHW} provided some spectral conditions for graphs to contain $\{K_2\}$-factors. Li and Miao \cite{LM} established
a lower bound on the adjacency spectral radius for a connected graph which ensures that this graph has a $P_{\geq2}$-factor. Zhou, Zhang and Sun \cite{ZZS}
presented an $A_{\alpha}$-spectral radius condition for a connected graph to contain a $P_{\geq2}$-factor. Zhou, Sun and Liu \cite{ZSL2} studied the
existence of a $P_{\geq2}$-factor in a connected graph and characterized a $P_{\geq2}$-factor in a connected graph via the distance signless Laplacian
spectral radius. Miao and Li \cite{ML} determined a lower bound on the adjacency spectral radius of a connected graph $G$ to guarantee that $G$ has a
$\{K_{1,j}|1\leq j\leq k\}$-factor, and presented an upper bound on the distance spectral radius of a connected graph $G$ to ensure that $G$ contains a
$\{K_{1,j}|1\leq j\leq k\}$-factor.

Motivated by \cite{O,KLY} directly, we investigate the existence of $\{P_3,P_4,P_5\}$-factors in connected graphs, and establish a relationship between
$A_{\alpha}$-spectral radius and $\{P_3,P_4,P_5\}$-factors in connected graphs. Our main result is shown as follows.

\medskip

\noindent{\textbf{Theorem 1.1.}} Let $\alpha$ be a real number with $0\leq\alpha<\frac{2}{3}$, and let $G$ be a connected graph of order $n$ with $n\geq25$.
If $G$ satisfies
$$
\lambda_{\alpha}(G)\geq\lambda_{\alpha}(K_1\vee(K_{n-2}\cup K_1)),
$$
then $G$ has a $\{P_3,P_4,P_5\}$-factor unless $G=K_1\vee(K_{n-2}\cup K_1)$.

\medskip

In fact, a $\{P_3,P_4,P_5\}$-factor is also a $P_{\geq3}$-factor. Then the following corollary holds.

\medskip

\noindent{\textbf{Corollary 1.2.}} Let $\alpha$ be a real number with $0\leq\alpha<\frac{2}{3}$, and let $G$ be a connected graph of order $n$ with $n\geq25$.
If $G$ satisfies
$$
\lambda_{\alpha}(G)\geq\lambda_{\alpha}(K_1\vee(K_{n-2}\cup K_1)),
$$
then $G$ has a $P_{\geq3}$-factor unless $G=K_1\vee(K_{n-2}\cup K_1)$.

\section{Some preliminaries}

In 2010, Kano, Lu and Yu \cite{KLY} provided a sufficient condition for the existence of $\{P_3,P_4,P_5\}$-factors in graphs.

\medskip

\noindent{\textbf{Lemma 2.1}} (Kano, Lu and Yu \cite{KLY}). If a graph $G$ satisfies
$$
i(G-S)\leq\frac{2}{3}|S|
$$
for any subset $S\subset V(G)$, then $G$ contains a $\{P_3,P_4,P_5\}$-factor.

\medskip

\noindent{\textbf{Lemma 2.2}} (Nikiforov \cite{N}). For a complete graph $K_n$, we conclude
$$
\lambda_{\alpha}(K_n)=n-1.
$$

\medskip

\noindent{\textbf{Lemma 2.3}} (Nikiforov \cite{N}). If $G$ is a connected graph, and $H$ is a proper subgraph of $G$, then we have
$$
\lambda_{\alpha}(G)>\lambda_{\alpha}(H).
$$

\medskip

Let $M$ be a real symmetric matrix of order $n$ whose columns and rows are indexed by $V=\{1,2,\ldots,n\}$, where $V=V_1\cup V_2\cup\cdots\cup V_t$,
$|V_i|=n_i$ and $n=\sum\limits_{i=1}^{t}n_i$. Assume that $M$ is a matrix with the partition $\pi: V=V_1\cup V_2\cup\cdots\cup V_t$, that is,
\begin{align*}
M=\left(
  \begin{array}{cccc}
    M_{11} & M_{12} & \cdots & M_{1t}\\
    M_{21} & M_{22} & \cdots & M_{2t}\\
    \vdots & \vdots & \ddots & \vdots\\
    M_{t1} & M_{t2} &\cdots & M_{tt}\\
  \end{array}
\right),
\end{align*}
where $M_{ij}$ denotes the submatrix (block) of $M$ formed by rows in $V_i$ and columns in $V_j$. The average row sum of $M_{ij}$ is denoted by $m_{ij}$.
Then the matrix $M_{\pi}=(m_{ij})$ is called the quotient matrix of $M$. In particular, if the row sum of each block $M_{ij}$ is a constant, then the
partition is called equitable.

\medskip

\noindent{\textbf{Lemma 2.4}} (You, Yang, So and Xi \cite{YYSX}). Let $M$ be a real matrix with an equitable partition $\pi$, and let $M_{\pi}$ be the
corresponding quotient matrix. Then every eigenvalue of $M_{\pi}$ is an eigenvalue of $M$. Furthermore, if $M$ is a nonnegative matrix, then the largest
eigenvalues of $M$ is equal to the largest eigenvalues of $M_{\pi}$.

\medskip

\noindent{\textbf{Lemma 2.5}} (Haemers \cite{Hi}). Let $M$ be a Hermitian matrix of order $s$, and let $N$ be a principal submatrix of $M$ of order $t$.
If $\lambda_1\geq\lambda_2\geq\cdots\geq\lambda_s$ are the eigenvalues of $M$ and $\mu_1\geq\mu_2\geq\cdots\geq\mu_t$ are the eigenvalues of $N$, then
$\lambda_i\geq\mu_i\geq\lambda_{s-t+i}$ for $1\leq i\leq t$.

\section{The proof of Theorem 1.1}

\noindent{\it Proof of Theorem 1.1.} Suppose, to the contrary, that $G$ contains no $\{P_3,P_4,P_5\}$-factor. Then it follows from Lemma 2.1 that
$i(G-S)>\frac{2}{3}|S|$ for some nonempty subset $S\subset V(G)$. Let $|S|=s$ and $i(G-S)=i$. According to the integrity of $i(G-S)$, we get
$i\geq\lfloor\frac{2}{3}s\rfloor+1$. From the above discussion, we easily see that $G$ is a spanning subgraph of
$G_1=K_s\vee(K_{n-\lfloor\frac{5}{3}s\rfloor-1}\cup(\lfloor\frac{2}{3}s\rfloor+1)K_1)$. Combining this with Lemma 2.3, we deduce
\begin{align}\label{eq:3.1}
\lambda_{\alpha}(G)\leq\lambda_{\alpha}(G_1),
\end{align}
with equality holding if and only if $G=G_1$. The following proof will be divided into three cases by the value of $n$.

\noindent{\bf Case 1.} $n\geq\lfloor\frac{5}{3}s\rfloor+3$.

Recall that $G_1=K_s\vee(K_{n-\lfloor\frac{5}{3}s\rfloor-1}\cup(\lfloor\frac{2}{3}s\rfloor+1)K_1)$. The quotient matrix of $A_{\alpha}(G_1)$ by the
partition $V(G_1)=V(K_s)\cup V(K_{n-\lfloor\frac{5}{3}s\rfloor-1})\cup V((\lfloor\frac{2}{3}s\rfloor+1)K_1)$ can be written as
\begin{align*}
B_1=\left(
  \begin{array}{ccc}
    \alpha n-\alpha s+s-1 & (1-\alpha)(n-\lfloor\frac{5}{3}s\rfloor-1) & (1-\alpha)(\lfloor\frac{2}{3}s\rfloor+1)\\
    (1-\alpha)s & n+\alpha s-\lfloor\frac{5}{3}s\rfloor-2 & 0\\
    (1-\alpha)s & 0 & \alpha s\\
  \end{array}
\right).
\end{align*}
By a direct calculation, the characteristic polynomial of $B_1$ is
\begin{align}\label{eq:3.2}
\varphi_{B_1}(x)=&x^{3}-\Big(\alpha n+n+\alpha s-\Big\lfloor\frac{2}{3}s\Big\rfloor-3\Big)x^{2}\nonumber\\
&+\Big(\alpha n^{2}+\alpha^{2}sn+\alpha sn-\alpha n\Big\lfloor\frac{2}{3}s\Big\rfloor-2\alpha n-n-s\Big\lfloor\frac{2}{3}s\Big\rfloor-2\alpha s-s+\Big\lfloor\frac{2}{3}s\Big\rfloor+2\Big)x\nonumber\\
&-\alpha^{2}sn^{2}+2\alpha^{2}sn\Big\lfloor\frac{2}{3}s\Big\rfloor-2\alpha sn\Big\lfloor\frac{2}{3}s\Big\rfloor+3\alpha^{2}sn-\alpha sn+sn\Big\lfloor\frac{2}{3}s\Big\rfloor+sn\nonumber\\
&-2\alpha^{2}s^{2}\Big\lfloor\frac{2}{3}s\Big\rfloor-\alpha^{2}s\Big\lfloor\frac{2}{3}s\Big\rfloor^{2}+3\alpha s^{2}\Big\lfloor\frac{2}{3}s\Big\rfloor+2\alpha s\Big\lfloor\frac{2}{3}s\Big\rfloor^{2}
-s^{2}\Big\lfloor\frac{2}{3}s\Big\rfloor-s\Big\lfloor\frac{2}{3}s\Big\rfloor^{2}-2\alpha^{2}s^{2}\nonumber\\
&-3\alpha^{2}s\Big\lfloor\frac{2}{3}s\Big\rfloor+3\alpha s^{2}+5\alpha s\Big\lfloor\frac{2}{3}s\Big\rfloor-s^{2}-3s\Big\lfloor\frac{2}{3}s\Big\rfloor
-2\alpha^{2}s+2\alpha s-2s.
\end{align}
Notice that the partition $V(G_1)=V(K_s)\cup V(K_{n-\lfloor\frac{5}{3}s\rfloor-1})\cup V((\lfloor\frac{2}{3}s\rfloor+1)K_1)$ is equitable. By virtue of
Lemma 2.4, $\lambda_{\alpha}(G_1)$ is the largest root of $\varphi_{B_1}(x)=0$, that is, $\varphi_{B_1}(\lambda_{\alpha}(G_1))=0$. Let $\theta_1=\lambda_{\alpha}(G_1)\geq\theta_2\geq\theta_3$ be the three roots of $\varphi_{B_1}(x)=0$ and
$Q=\mbox{diag}(s,n-\lfloor\frac{5}{3}s\rfloor-1,\lfloor\frac{2}{3}s\rfloor+1)$. We easily see that
$$
Q^{\frac{1}{2}}B_1Q^{-\frac{1}{2}}=\left(
  \begin{array}{ccc}
    \alpha n-\alpha s+s-1 & (1-\alpha)s^{\frac{1}{2}}(n-\lfloor\frac{5}{3}s\rfloor-1)^{\frac{1}{2}} & (1-\alpha)s^{\frac{1}{2}}(\lfloor\frac{2}{3}s\rfloor+1)^{\frac{1}{2}}\\
    (1-\alpha)s^{\frac{1}{2}}(n-\lfloor\frac{5}{3}s\rfloor-1)^{\frac{1}{2}} & n+\alpha s-\lfloor\frac{5}{3}s\rfloor-2 & 0\\
    (1-\alpha)s^{\frac{1}{2}}(\lfloor\frac{2}{3}s\rfloor+1)^{\frac{1}{2}} & 0 & \alpha s\\
  \end{array}
\right)
$$
is symmetric, and
\begin{align*}
\left(
  \begin{array}{ccc}
    n+\alpha s-\lfloor\frac{5}{3}s\rfloor-2 & 0\\
    0 & \alpha s\\
  \end{array}
\right)
\end{align*}
is a submatrix of $Q^{\frac{1}{2}}B_1Q^{-\frac{1}{2}}$. Since $Q^{\frac{1}{2}}B_1Q^{-\frac{1}{2}}$ and $B_1$ have the same eigenvalues, Lemma 2.5 (the
Cauchy interlacing theorem) leads to
\begin{align}\label{eq:3.3}
\theta_2\leq n+\alpha s-\Big\lfloor\frac{5}{3}s\Big\rfloor-2<\left\{
\begin{array}{ll}
n-4,&\mbox{if} \ s\equiv0 \ (\mbox{mod} \ 3);\\
n-2,&\mbox{if} \ s\equiv1 \ (\mbox{mod} \ 3);\\
n-3,&\mbox{if} \ s\equiv2 \ (\mbox{mod} \ 3).\\
\end{array}
\right.
\end{align}

For the graph $G_*=K_1\vee(K_{n-2}\cup K_1)$, its adjacency matrix $A(G_*)$ admits the quotient matrix $B_*$ which is derived by replacing $s$ with 1
in $B_1$, and $B_*$ admits the characteristic polynomial $\varphi_{B_*}(x)$ which is derived by replacing $s$ with 1 in $\varphi_{B_1}(x)$. Hence, we
have
\begin{align*}
\varphi_{B_*}(x)=&x^{3}-(\alpha n+n+\alpha-3)x^{2}+(\alpha n^{2}+\alpha^{2}n-\alpha n-n-2\alpha+1)x\\
&-\alpha^{2}n^{2}+3\alpha^{2}n-\alpha n+n-4\alpha^{2}+5\alpha-3.
\end{align*}
In view of Lemma 2.4, $\lambda_{\alpha}(G_*)$ is the largest root of $\varphi_{B_*}(x)=0$, that is, $\varphi_{B_*}(\lambda_{\alpha}(G_*))=0$. If
$s=1$, then $G_1=G_*$, and so $\lambda_{\alpha}(G_1)=\lambda_{\alpha}(G_*)=\lambda_{\alpha}(K_1\vee(K_{n-2}\cup K_1))$. Combining this with \eqref{eq:3.1},
we deduce $\lambda_{\alpha}(G)\leq\lambda_{\alpha}(K_1\vee(K_{n-2}\cup K_1))$, where the equality holds if and only if $G=K_1\vee(K_{n-2}\cup K_1)$.
This is a contradiction. Next, we consider $s\geq2$.

Since $K_{n-1}$ is a proper subgraph of $G_*=K_1\vee(K_{n-2}\cup K_1)$, it follows from \eqref{eq:3.3}, Lemmas 2.2 and 2.3 that
\begin{align}\label{eq:3.4}
\lambda_{\alpha}(K_1\vee(K_{n-2}\cup K_1))>\lambda_{\alpha}(K_{n-1})=n-2>\theta_2.
\end{align}
We shall consider three subcases by the value of $s$.

\noindent{\bf Subcase 1.1.} $s\equiv0$ (mod 3).

In this subcase, $s\geq3$, $\lfloor\frac{2}{3}s\rfloor=\frac{2}{3}s$ and $n\geq\lfloor\frac{5}{3}s\rfloor+3=\frac{5}{3}s+3$. According to \eqref{eq:3.2},
we admit
\begin{align*}
\varphi_{B_1}(x)=&x^{3}-\Big(\alpha n+n+\alpha s-\frac{2}{3}s-3\Big)x^{2}\\
&+\Big(\alpha n^{2}+\alpha^{2}sn+\frac{1}{3}\alpha sn-2\alpha n-n-\frac{2}{3}s^{2}-2\alpha s-\frac{1}{3}s+2\Big)x\\
&-\alpha^{2}sn^{2}+\frac{4}{3}\alpha^{2}s^{2}n-\frac{4}{3}\alpha s^{2}n+3\alpha^{2}sn-\alpha sn+\frac{2}{3}s^{2}n+sn-\frac{16}{9}\alpha^{2}s^{3}\\
&+\frac{26}{9}\alpha s^{3}-\frac{10}{9}s^{3}-4\alpha^{2}s^{2}+\frac{19}{3}\alpha s^{2}-3s^{2}-2\alpha^{2}s+2\alpha s-2s.
\end{align*}
Let $G_2=K_3\vee(K_{n-6}\cup3K_1)$. Then its adjacency matrix $A(G_2)$ has the quotient matrix $B_2$ which is derived by replacing $s$ with 3 in $B_1$,
and $B_2$ has the characteristic polynomial $\varphi_{B_2}(x)$ which is obtained by replacing $s$ with 3 in $\varphi_{B_1}(x)$. Therefore, we obtain
\begin{align*}
\varphi_{B_2}(x)=&x^{3}-(\alpha n+n+3\alpha-5)x^{2}+(\alpha n^{2}+3\alpha^{2}n-\alpha n-n-6\alpha-5)x\\
&-3\alpha^{2}n^{2}+21\alpha^{2}n-15\alpha n+9n-90\alpha^{2}+141\alpha-63.
\end{align*}
Using Lemma 2.4, $\lambda_{\alpha}(G_2)$ is the largest root of $\varphi_{B_2}(x)=0$, that is, $\varphi_{B_2}(\lambda_{\alpha}(G_2))=0$. We are to
verify $\lambda_{\alpha}(G_1)\leq\lambda_{\alpha}(G_2)$.

Since $K_{n-3}$ is a proper subgraph of $G_2=K_3\vee(K_{n-6}\cup3K_1)$, it follows from \eqref{eq:3.3}, Lemmas 2.2 and 2.3 that
\begin{align}\label{eq:3.5}
\lambda_{\alpha}(K_3\vee(K_{n-6}\cup3K_1))>\lambda_{\alpha}(K_{n-3})=n-4>\theta_2.
\end{align}
Write $\beta=\lambda_{\alpha}(K_3\vee(K_{n-6}\cup3K_1))$. Notice that $\varphi_{B_2}(\beta)=0$. By a direct computation, we get
\begin{align}\label{eq:3.6}
\varphi_{B_1}(\beta)=\varphi_{B_1}(\beta)-\varphi_{B_2}(\beta)=\frac{1}{9}(s-3)f_1(\beta),
\end{align}
where $f_1(\beta)=(6-9\alpha)\beta^{2}+(9\alpha^{2}n+3\alpha n-6s-18\alpha-21)\beta-9\alpha^{2}n^{2}+3\alpha^{2}n(4s+21)-3\alpha n(4s+15)+3n(2s+9)
-2\alpha^{2}(8s^{2}+42s+135)+\alpha(26s^{2}+135s+423)-10s^{2}-57s-189$. Notice that
\begin{align}\label{eq:3.7}
-\frac{9\alpha^{2}n+3\alpha n-6s-18\alpha-21}{2(6-9\alpha)}<n-4<\beta
\end{align}
by \eqref{eq:3.5}, $s\geq6$ and $n\geq\frac{5}{3}s+3$. Since the symmetry axis of $f_1(\beta)$ is $\beta=-\frac{9\alpha^{2}n+3\alpha n-6s-18\alpha-21}{2(6-9\alpha)}$,
it follows from \eqref{eq:3.7} that
\begin{align}\label{eq:3.8}
f_1(\beta)>&f_1(n-4)\nonumber\\
=&(6-6\alpha)n^{2}+(12\alpha^{2}s-12\alpha s+27\alpha^{2}-3\alpha-42)n\nonumber\\
&-2\alpha^{2}(8s^{2}+42s+135)+\alpha(26s^{2}+135s+351)-10s^{2}-33s-9.
\end{align}
Let $f_2(n)=(6-6\alpha)n^{2}+(12\alpha^{2}s-12\alpha s+27\alpha^{2}-3\alpha-42)n-2\alpha^{2}(8s^{2}+42s+135)+\alpha(26s^{2}+135s+351)-10s^{2}-33s-9$.
Note that
$$
-\frac{12\alpha^{2}s-12\alpha s+27\alpha^{2}-3\alpha-42}{2(6-6\alpha)}<\frac{5}{3}s+3\leq n
$$
by $s\geq15$ and $0\leq\alpha<\frac{2}{3}$. Thus, we deduce
\begin{align}\label{eq:3.9}
f_2(n)\geq&f_2\Big(\frac{5}{3}s+3\Big)\nonumber\\
=&\frac{1}{3}((12s^{2}-9s-567)\alpha^{2}+(-32s^{2}+102s+864)\alpha+20s^{2}-129s-243)\nonumber\\
>&\frac{1}{3}\Big(\frac{4}{9}(12s^{2}-9s-567)+\frac{2}{3}(-32s^{2}+102s+864)+20s^{2}-129s-243\Big)\nonumber\\
=&\frac{1}{9}(12s^{2}-195s+243)\nonumber\\
>&0,
\end{align}
where the last two inequalities hold from $\frac{32s^{2}-102s-864}{2(12s^{2}-9s-567)}>\frac{2}{3}>\alpha\geq0$ and $s\geq15$, respectively.

If $s\in\{6,9,12\}$, then
\begin{align*}
-\frac{12\alpha^{2}s-12\alpha s+27\alpha^{2}-3\alpha-42}{2(6-6\alpha)}=&\left\{
\begin{array}{ll}
\frac{14+25\alpha-33\alpha^{2}}{4-4\alpha},&\mbox{if} \ s=6\\
\frac{14+37\alpha-45\alpha^{2}}{4-4\alpha},&\mbox{if} \ s=9\\
\frac{14+49\alpha-57\alpha^{2}}{4-4\alpha},&\mbox{if} \ s=12\\
\end{array}
\right.\\
<&25\leq n
\end{align*}
by $0\leq\alpha<\frac{2}{3}$. Thus, we obtain
\begin{align}\label{eq:3.10}
f_2(n)\geq&f_2(25)\nonumber\\
=&(-16s^{2}+216s+405)\alpha^{2}+(26s^{2}-165s-3474)\alpha-10s^{2}-33s+2691\nonumber\\
=&\left\{
\begin{array}{ll}
1125\alpha^{2}-3528\alpha+2133,&\mbox{if} \ s=6\\
1053\alpha^{2}-2853\alpha+1584,&\mbox{if} \ s=9\\
693\alpha^{2}-1710\alpha+855,&\mbox{if} \ s=12\\
\end{array}
\right.\nonumber\\
>&0
\end{align}
by $0\leq\alpha<\frac{2}{3}$.

From \eqref{eq:3.9} and \eqref{eq:3.10}, we infer $f_2(n)>0$ for $s\geq6$ and $s\equiv0$ (mod 3). Combining this with \eqref{eq:3.6} and \eqref{eq:3.8},
we conclude
\begin{align}\label{eq:3.11}
\varphi_{B_1}(\beta)=\frac{1}{9}(s-3)f_1(\beta)\geq\frac{1}{9}(s-3)f_1(n-4)=\frac{1}{9}(s-3)f_2(n)\geq0
\end{align}
for $s\geq3$ and $s\equiv0$ (mod 3). Recall that $\lambda_{\alpha}(G_1)$ is the largest root of $\varphi_{B_1}(x)=0$. As $\theta_2<n-4<\lambda_{\alpha}(K_3\vee(K_{n-6}\cup3K_1))=\beta$ (see \eqref{eq:3.5}), we deduce
\begin{align}\label{eq:3.12}
\lambda_{\alpha}(G_1)\leq\beta=\lambda_{\alpha}(K_3\vee(K_{n-6}\cup3K_1))=\lambda_{\alpha}(G_2)
\end{align}
by \eqref{eq:3.11}.

In what follows, we are to show $\lambda_{\alpha}(G_2)<n-2$. By a direct calculation, we get
\begin{align*}
\varphi_{B_2}(n-2)=&(n-2)^{3}-(\alpha n+n+3\alpha-5)(n-2)^{2}\\
&+(\alpha n^{2}+3\alpha^{2}n-\alpha n-n-6\alpha-5)(n-2)\\
&-3\alpha^{2}n^{2}+21\alpha^{2}n-15\alpha n+9n-90\alpha^{2}+141\alpha-63\\
=&(2-2\alpha)n^{2}+(15\alpha^{2}-11\alpha-6)n-90\alpha^{2}+141\alpha-41\\
\geq&(2-2\alpha)(25)^{2}+25(15\alpha^{2}-11\alpha-6)-90\alpha^{2}+141\alpha-41\\
=&285\alpha^{2}-1384\alpha+1059\\
>&0,
\end{align*}
where the last two inequalities hold from $-\frac{15\alpha^{2}-11\alpha-6}{2(2-2\alpha)}<25\leq n$ and $\frac{1384}{2\times 285}>\frac{2}{3}>\alpha\geq0$,
respectively. Hence, we infer
\begin{align}\label{eq:3.13}
\lambda_{\alpha}(G_2)<n-2.
\end{align}

According to \eqref{eq:3.1}, \eqref{eq:3.4}, \eqref{eq:3.12} and \eqref{eq:3.13}, we have
$$
\lambda_{\alpha}(G)\leq\lambda_{\alpha}(G_1)\leq\lambda_{\alpha}(G_2)<n-2<\lambda_{\alpha}(K_1\vee(K_{n-2}\cup K_1)),
$$
which contradicts $\lambda_{\alpha}(G)\geq\lambda_{\alpha}(K_1\vee(K_{n-2}\cup K_1))$.

\noindent{\bf Subcase 1.2.} $s\equiv1$ (mod 3).

In this subcase, $s\geq4$, $\lfloor\frac{2}{3}s\rfloor=\frac{2s-2}{3}$ and $n\geq\lfloor\frac{5}{3}s\rfloor+3=\frac{5s+7}{3}$. In view of \eqref{eq:3.2},
we obtain
\begin{align*}
\varphi_{B_1}(x)=&x^{3}-\Big(\alpha n+n+\alpha s-\frac{2}{3}s-\frac{7}{3}\Big)x^{2}\\
&+\Big(\alpha n^{2}+\alpha^{2}sn+\frac{1}{3}\alpha sn-\frac{4}{3}\alpha n-n-\frac{2}{3}s^{2}-2\alpha s+\frac{1}{3}s+\frac{4}{3}\Big)x\\
&-\alpha^{2}sn^{2}+\frac{4}{3}\alpha^{2}s^{2}n-\frac{4}{3}\alpha s^{2}n+\frac{5}{3}\alpha^{2}sn+\frac{1}{3}\alpha sn+\frac{2}{3}s^{2}n+\frac{1}{3}sn-\frac{16}{9}\alpha^{2}s^{3}\\
&+\frac{26}{9}\alpha s^{3}-\frac{10}{9}s^{3}-\frac{16}{9}\alpha^{2}s^{2}+\frac{23}{9}\alpha s^{2}-\frac{13}{9}s^{2}-\frac{4}{9}\alpha^{2}s-\frac{4}{9}\alpha s-\frac{4}{9}s.
\end{align*}
Write $\gamma=\lambda_{\alpha}(K_1\vee(K_{n-2}\cup K_1))$. Notice that $\varphi_{B_*}(\gamma)=0$. A simple calculation yields that
\begin{align}\label{eq:3.14}
\varphi_{B_1}(\gamma)=\varphi_{B_1}(\gamma)-\varphi_{B_*}(\gamma)=\frac{1}{9}(s-1)g_1(\gamma),
\end{align}
where $g_1(\gamma)=(6-9\alpha)\gamma^{2}+(9\alpha^{2}n+3\alpha n-6s-18\alpha-3)\gamma-9\alpha^{2}n^{2}+3\alpha^{2}n(4s+9)-3\alpha n(4s+3)+3n(2s+3)
-4\alpha^{2}(4s^{2}+8s+9)+\alpha(26s^{2}+49s+45)-10s^{2}-23s-27$. Note that
$$
-\frac{9\alpha^{2}n+3\alpha n-6s-18\alpha-3}{2(6-9\alpha)}<n-2<\gamma
$$
by \eqref{eq:3.4}, $s\geq4$ and $n\geq\frac{5s+7}{3}$. Thus, we have
\begin{align}\label{eq:3.15}
g_1(\gamma)>&g_1(n-2)\nonumber\\
=&(6-6\alpha)n^{2}+(12\alpha^{2}s-12\alpha s+9\alpha^{2}+3\alpha-18)n\nonumber\\
&-4\alpha^{2}(4s^{2}+8s+9)+\alpha(26s^{2}+49s+45)-10s^{2}-11s+3.
\end{align}
Let $g_2(n)=(6-6\alpha)n^{2}+(12\alpha^{2}s-12\alpha s+9\alpha^{2}+3\alpha-18)n-4\alpha^{2}(4s^{2}+8s+9)+\alpha(26s^{2}+49s+45)-10s^{2}-11s+3$. It follows
from $s\geq4$ and $0\leq\alpha<\frac{2}{3}$ that
$$
-\frac{12\alpha^{2}s-12\alpha s+9\alpha^{2}+3\alpha-18}{2(6-6\alpha)}<\frac{5s+7}{3}\leq n,
$$
and so
\begin{align}\label{eq:3.16}
g_2(n)\geq&g_2\Big(\frac{5s+7}{3}\Big)\nonumber\\
=&\frac{1}{3}((12s^{2}+33s-45)\alpha^{2}-(32s^{2}+62s-58)\alpha+20s^{2}+17s-19)\nonumber\\
>&\frac{1}{3}\Big(\frac{4}{9}(12s^{2}+33s-45)-\frac{2}{3}(32s^{2}+62s-58)+20s^{2}+17s-19\Big)\nonumber\\
=&\frac{1}{9}(12s^{2}-29s-1)\nonumber\\
>&0,
\end{align}
where the last two inequalities hold from $\frac{32s^{2}+62s-58}{2(12s^{2}+33s-45)}>\frac{2}{3}>\alpha\geq0$ and $s\geq4$, respectively.

By virtue of \eqref{eq:3.14}, \eqref{eq:3.15} and \eqref{eq:3.16}, we obtain
\begin{align}\label{eq:3.17}
\varphi_{B_1}(\gamma)=\frac{1}{9}(s-1)g_1(\gamma)>\frac{1}{9}(s-1)g_1(n-2)=\frac{1}{9}(s-1)g_2(n)>0
\end{align}
for $s\geq4$ and $s\equiv1$ (mod 3). Notice that $\lambda_{\alpha}(G_1)$ is the largest root of $\varphi_{B_1}(x)=0$. As $\theta_2<n-2<\lambda_{\alpha}(K_1\vee(K_{n-2}\cup K_1))=\gamma$ (see \eqref{eq:3.4}), we conclude
$$
\lambda_{\alpha}(G_1)<\gamma=\lambda_{\alpha}(K_1\vee(K_{n-2}\cup K_1))
$$
by \eqref{eq:3.17}. Combining this with \eqref{eq:3.1}, we have
$$
\lambda_{\alpha}(G)\leq\lambda_{\alpha}(G_1)<\lambda_{\alpha}(K_1\vee(K_{n-2}\cup K_1)),
$$
which is a contradiction to $\lambda_{\alpha}(G)\geq\lambda_{\alpha}(K_1\vee(K_{n-2}\cup K_1))$.

\noindent{\bf Subcase 1.3.} $s\equiv2$ (mod 3).

In this subcase, $s\geq2$, $\lfloor\frac{2}{3}s\rfloor=\frac{2s-1}{3}$ and $n\geq\lfloor\frac{5}{3}s\rfloor+3=\frac{5s+8}{3}$. Using \eqref{eq:3.2}, we
possess
\begin{align*}
\varphi_{B_1}(x)=&x^{3}-\Big(\alpha n+n+\alpha s-\frac{2}{3}s-\frac{8}{3}\Big)x^{2}\\
&+\Big(\alpha n^{2}+\alpha^{2}sn+\frac{1}{3}\alpha sn-\frac{5}{3}\alpha n-n-\frac{2}{3}s^{2}-2\alpha s+\frac{5}{3}\Big)x\\
&-\alpha^{2}sn^{2}+\frac{4}{3}\alpha^{2}s^{2}n-\frac{4}{3}\alpha s^{2}n+\frac{7}{3}\alpha^{2}sn-\frac{1}{3}\alpha sn+\frac{2}{3}s^{2}n+\frac{2}{3}sn-\frac{16}{9}\alpha^{2}s^{3}\\
&+\frac{26}{9}\alpha s^{3}-\frac{10}{9}s^{3}-\frac{26}{9}\alpha^{2}s^{2}+\frac{40}{9}\alpha s^{2}-\frac{20}{9}s^{2}-\frac{10}{9}\alpha^{2}s+\frac{5}{9}\alpha s-\frac{10}{9}s.
\end{align*}
Let $G_3=K_2\vee(K_{n-4}\cup2K_1)$. Then its adjacency matrix $A(G_3)$ has the quotient matrix $B_3$ which is obtained by replacing $s$ with 2 in $B_1$,
and $B_3$ admits the characteristic polynomial $\varphi_{B_3}(x)$ which is derived by replacing $s$ with 2 in $\varphi_{B_1}(x)$. Hence, we get
\begin{align*}
\varphi_{B_3}(x)=&x^{3}-(\alpha n+n+2\alpha-4)x^{2}+(\alpha n^{2}+2\alpha^{2}n-\alpha n-n-4\alpha-1)x\\
&-2\alpha^{2}n^{2}+10\alpha^{2}n-6\alpha n+4n-28\alpha^{2}+42\alpha-20.
\end{align*}
In terms of Lemma 2.4, $\lambda_{\alpha}(G_3)$ is the largest root of $\varphi_{B_3}(x)=0$, that is, $\varphi_{B_3}(\lambda_{\alpha}(G_3))=0$. We are to
verify $\lambda_{\alpha}(G_1)\leq\lambda_{\alpha}(G_3)$.

Note that $K_{n-2}$ is a proper subgraph of $G_3=K_2\vee(K_{n-4}\cup2K_1)$. By means of \eqref{eq:3.3}, Lemmas 2.2 and 2.3, we conclude
\begin{align}\label{eq:3.18}
\lambda_{\alpha}(K_2\vee(K_{n-4}\cup2K_1))>\lambda_{\alpha}(K_{n-2})=n-3>\theta_2.
\end{align}
Write $\eta=\lambda_{\alpha}(K_2\vee(K_{n-4}\cup2K_1))$. Note that $\varphi_{B_3}(\eta)=0$. By a direct computation, we possess
\begin{align}\label{eq:3.19}
\varphi_{B_1}(\eta)=\varphi_{B_1}(\eta)-\varphi_{B_3}(\eta)=\frac{1}{9}(s-2)h_1(\eta),
\end{align}
where $h_1(\eta)=(6-9\alpha)\eta^{2}+(9\alpha^{2}n+3\alpha n-6s-18\alpha-12)\eta-9\alpha^{2}n^{2}+3\alpha^{2}n(4s+15)-3\alpha n(4s+9)+6n(s+3)
-2\alpha^{2}(8s^{2}+29s+63)+\alpha(26s^{2}+92s+189)-10s^{2}-40s-90$. According to \eqref{eq:3.18}, $s\geq5$ and $n\geq\frac{5s+8}{3}$, we deduce
$$
-\frac{9\alpha^{2}n+3\alpha n-6s-18\alpha-12}{2(6-9\alpha)}<n-3<\eta,
$$
and so
\begin{align}\label{eq:3.20}
h_1(\eta)>&h_1(n-3)\nonumber\\
=&(6-6\alpha)n^{2}+(12\alpha^{2}s-12\alpha s+18\alpha^{2}-30)n\nonumber\\
&-2\alpha^{2}(8s^{2}+29s+63)+\alpha(26s^{2}+92s+162)-10s^{2}-22s.
\end{align}
Let $h_2(n)=(6-6\alpha)n^{2}+(12\alpha^{2}s-12\alpha s+18\alpha^{2}-30)n-2\alpha^{2}(8s^{2}+29s+63)+\alpha(26s^{2}+92s+162)-10s^{2}-22s$. It follows
from $s\geq11$ and $0\leq\alpha<\frac{2}{3}$ that
$$
-\frac{12\alpha^{2}s-12\alpha s+18\alpha^{2}-30}{2(6-6\alpha)}<\frac{5s+8}{3}\leq n,
$$
and so
\begin{align}\label{eq:3.21}
h_2(n)\geq&h_2\Big(\frac{5s+8}{3}\Big)\nonumber\\
=&\frac{1}{3}((12s^{2}+12s-234)\alpha^{2}+(-32s^{2}+20s+358)\alpha+20s^{2}-56s-112)\nonumber\\
>&\frac{1}{3}\Big(\frac{4}{9}(12s^{2}+12s-234)+\frac{2}{3}(-32s^{2}+20s+358)+20s^{2}-56s-112\Big)\nonumber\\
=&\frac{1}{9}(12s^{2}-112s+68)\nonumber\\
>&0,
\end{align}
where the last two inequalities hold from $\frac{32s^{2}-20s-358}{2(12s^{2}+12s-234)}>\frac{2}{3}>\alpha\geq0$ and $s\geq11$, respectively.

If $s\in\{5,8\}$, then
\begin{align*}
-\frac{12\alpha^{2}s-12\alpha s+18\alpha^{2}-30}{2(6-6\alpha)}=&\left\{
\begin{array}{ll}
\frac{5+10\alpha-13\alpha^{2}}{2-2\alpha},&\mbox{if} \ s=5\\
\frac{5+16\alpha-19\alpha^{2}}{2-2\alpha},&\mbox{if} \ s=8\\
\end{array}
\right.\\
<&25\leq n
\end{align*}
by $0\leq\alpha<\frac{2}{3}$. Hence, we infer
\begin{align}\label{eq:3.22}
h_2(n)\geq&h_2(25)\nonumber\\
=&(-16s^{2}+242s+324)\alpha^{2}+(26s^{2}-208s-3588)\alpha-10s^{2}-22s+3000\nonumber\\
=&\left\{
\begin{array}{ll}
1134\alpha^{2}-3978\alpha+2640,&\mbox{if} \ s=5\\
1236\alpha^{2}-3588\alpha+2184,&\mbox{if} \ s=8\\
\end{array}
\right.\nonumber\\
>&0
\end{align}
by $0\leq\alpha<\frac{2}{3}$.

According to \eqref{eq:3.21} and \eqref{eq:3.22}, we conclude $h_2(n)>0$ for $s\geq5$ and $s\equiv2$ (mod 3). Together with \eqref{eq:3.19} and \eqref{eq:3.20},
we get
\begin{align}\label{eq:3.23}
\varphi_{B_1}(\eta)=\frac{1}{9}(s-2)h_1(\eta)\geq\frac{1}{9}(s-2)h_1(n-3)=\frac{1}{9}(s-2)h_2(n)\geq0
\end{align}
for $s\geq2$ and $s\equiv2$ (mod 3). Recall that $\lambda_{\alpha}(G_1)$ is the largest root of $\varphi_{B_1}(x)=0$. As $\theta_2<n-3<\lambda_{\alpha}(K_2\vee(K_{n-4}\cup2K_1))=\eta$ (see \eqref{eq:3.18}), we obtain
\begin{align}\label{eq:3.24}
\lambda_{\alpha}(G_1)\leq\eta=\lambda_{\alpha}(K_2\vee(K_{n-4}\cup2K_1))=\lambda_{\alpha}(G_3)
\end{align}
by \eqref{eq:3.23}.

Next, we prove $\lambda_{\alpha}(G_3)<n-2$. A direct computation yields that
\begin{align*}
\varphi_{B_3}(n-2)=&(n-2)^{3}-(\alpha n+n+2\alpha-4)(n-2)^{2}\\
&+(\alpha n^{2}+2\alpha^{2}n-\alpha n-n-4\alpha-1)(n-2)\\
&-2\alpha^{2}n^{2}+10\alpha^{2}n-6\alpha n+4n-28\alpha^{2}+42\alpha-20\\
=&(1-\alpha)n^{2}+(6\alpha^{2}-4\alpha-3)n-28\alpha^{2}+42\alpha-10\\
\geq&(1-\alpha)(25)^{2}+25(6\alpha^{2}-4\alpha-3)-28\alpha^{2}+42\alpha-10\\
=&122\alpha^{2}-683\alpha+540\\
>&0,
\end{align*}
where the last two inequalities hold from $-\frac{6\alpha^{2}-4\alpha-3}{2(1-\alpha)}<25\leq n$ and $\frac{683}{2\times 122}>\frac{2}{3}>\alpha\geq0$,
respectively. Consequently, we deduce
\begin{align}\label{eq:3.25}
\lambda_{\alpha}(G_3)<n-2.
\end{align}

It follows from \eqref{eq:3.1}, \eqref{eq:3.4}, \eqref{eq:3.24} and \eqref{eq:3.25} that
$$
\lambda_{\alpha}(G)\leq\lambda_{\alpha}(G_1)\leq\lambda_{\alpha}(G_3)<n-2<\lambda_{\alpha}(K_1\vee(K_{n-2}\cup K_1)),
$$
which contradicts $\lambda_{\alpha}(G)\geq\lambda_{\alpha}(K_1\vee(K_{n-2}\cup K_1))$.

\noindent{\bf Case 2.} $n=\lfloor\frac{5}{3}s\rfloor+2$.

In this case, $G_1=K_s\vee(\lfloor\frac{2}{3}s\rfloor+2)K_1$. The quotient matrix of $A(G_1)$ with respect to the partition
$V(G_1)=V(K_s)\cup V((\lfloor\frac{2}{3}s\rfloor+2)K_1)$ equals
\begin{align*}
B_4=\left(
  \begin{array}{ccc}
    \alpha n-\alpha s+s-1 & (1-\alpha)(\lfloor\frac{2}{3}s\rfloor+2)\\
    (1-\alpha)s & \alpha s\\
  \end{array}
\right),
\end{align*}
for which we calculate the characteristic polynomial
$$
\varphi_{B_4}(x)=x^{2}-(\alpha n+s-1)x+\alpha^{2}sn-\alpha^{2}s^{2}+\alpha s^{2}-(1-\alpha)^{2}s\Big\lfloor\frac{2}{3}s\Big\rfloor-2\alpha^{2}s+3\alpha s-2s.
$$
Since the partition $V(G_1)=V(K_s)\cup V((\lfloor\frac{2}{3}s\rfloor+2)K_1)$ is equitable, it follows from Lemma 2.4 that $\lambda_{\alpha}(G_1)$ is
the largest root of $\varphi_{B_4}(x)=0$. Hence, we conclude
\begin{align}\label{eq:3.26}
\lambda_{\alpha}(G_1)=M,
\end{align}
where $M=\frac{\alpha n+s-1+\sqrt{(\alpha n+s-1)^{2}-4(\alpha^{2}sn-\alpha^{2}s^{2}+\alpha s^{2}-(1-\alpha)^{2}s\lfloor\frac{2}{3}s\rfloor-2\alpha^{2}s+3\alpha s-2s)}}{2}$. We are to prove $\lambda_{\alpha}(G_1)<n-2$. According to $n=\lfloor\frac{5}{3}s\rfloor+2$, we have
\begin{align}\label{eq:3.27}
(2(n-2)&-\alpha n-s+1)^{2}-(\alpha n+s-1)^{2}+4(\alpha^{2}sn-\alpha^{2}s^{2}+\alpha s^{2}-(1-\alpha)^{2}s\Big\lfloor\frac{2}{3}s\Big\rfloor-2\alpha^{2}s+3\alpha s-2s)\nonumber\\
=&(4-4\alpha)n^{2}+(4\alpha^{2}s-4s+8\alpha-12)n-4\alpha^{2}s^{2}+4\alpha s^{2}-4(1-\alpha)^{2}s\Big\lfloor\frac{2}{3}s\Big\rfloor-8\alpha^{2}s+12\alpha s+8\nonumber\\
=&\left\{
\begin{array}{ll}
\frac{4}{9}((4-4\alpha)s^{2}-(3\alpha+3)s),&\mbox{if} \ s\equiv0 \ (\mbox{mod 3});\\
\frac{4}{9}((4-4\alpha)s^{2}+(5\alpha-11)s+8\alpha-2),&\mbox{if} \ s\equiv1 \ (\mbox{mod 3});\\
\frac{4}{9}((4-4\alpha)s^{2}+(\alpha-7)s+5\alpha-2),&\mbox{if} \ s\equiv2 \ (\mbox{mod 3}).\\
\end{array}
\right.
\end{align}

\noindent{\bf Subcase 2.1.} $s\equiv0$ (mod 3).

Obviously, $n=\frac{5}{3}s+2\geq25$. Then $s\geq15$. Let $\psi_1(s)=(4-4\alpha)s^{2}-(3\alpha+3)s$. Note that $\frac{3\alpha+3}{2(4-4\alpha)}<15\leq s$.
Hence, we deduce
\begin{align}\label{eq:3.28}
\psi_1(s)\geq\psi_1(15)=9(95-105\alpha)>0.
\end{align}

\noindent{\bf Subcase 2.2.} $s\equiv1$ (mod 3).

It is obvious that $n=\frac{5s+4}{3}\geq25$. Then $s\geq16$. Let $\psi_2(s)=(4-4\alpha)s^{2}+(5\alpha-11)s+8\alpha-2$. Since $0\leq\alpha<\frac{2}{3}$
and $-\frac{5\alpha-11}{2(4-4\alpha)}<16\leq s$, we obtain
\begin{align}\label{eq:3.29}
\psi_2(s)\geq\psi_2(16)=18(47-52\alpha)>0.
\end{align}

\noindent{\bf Subcase 2.3.} $s\equiv2$ (mod 3).

Clearly, $n=\frac{5s+5}{3}\geq25$. Then $s\geq14$. Let $\psi_3(s)=(4-4\alpha)s^{2}+(\alpha-7)s+5\alpha-2$. Since $0\leq\alpha<\frac{2}{3}$ and
$-\frac{\alpha-7}{2(4-4\alpha)}<14\leq s$, we get
\begin{align}\label{eq:3.30}
\psi_3(s)\geq\psi_3(14)=9(76-85\alpha)>0.
\end{align}

According to \eqref{eq:3.26}, \eqref{eq:3.27}, \eqref{eq:3.28}, \eqref{eq:3.29} and \eqref{eq:3.30}, we conclude $\lambda_{\alpha}(G_1)<n-2$. Combining
this with \eqref{eq:3.1} and \eqref{eq:3.4}, we have $\lambda_{\alpha}(G)\leq\lambda_{\alpha}(G_1)<n-2<\lambda_{\alpha}(K_1\vee(K_{n-2}\cup K_1))$, which
contradicts $\lambda_{\alpha}(G)\geq\lambda_{\alpha}(K_1\vee(K_{n-2}\cup K_1))$.

\noindent{\bf Case 3.} $n=\lfloor\frac{5}{3}s\rfloor+1$.

In this case, $G_1=K_s\vee(\lfloor\frac{2}{3}s\rfloor+1)K_1$. Consider the partition $V(G_1)=V(K_s)\cup V((\lfloor\frac{2}{3}s\rfloor+1)K_1)$. The
corresponding quotient matrix of $A(G_1)$ equals
\begin{align*}
B_5=\left(
  \begin{array}{ccc}
    \alpha n-\alpha s+s-1 & (1-\alpha)(\lfloor\frac{2}{3}s\rfloor+1)\\
    (1-\alpha)s & \alpha s\\
  \end{array}
\right).
\end{align*}
Then the characteristic polynomial of $B_5$ is
$$
\varphi_{B_5}(x)=x^{2}-(\alpha n+s-1)x+\alpha^{2}sn-\alpha^{2}s^{2}+\alpha s^{2}-(1-\alpha)^{2}s\Big\lfloor\frac{2}{3}s\Big\rfloor-\alpha^{2}s+\alpha s-s.
$$
Since the partition $V(G_1)=V(K_s)\cup V((\lfloor\frac{2}{3}s\rfloor+1)K_1)$ is equitable, $\lambda_{\alpha}(G_1)$ is the largest root of
$\varphi_{B_5}(x)=0$ by Lemma 2.4. Thus, we obtain
\begin{align}\label{eq:3.31}
\lambda_{\alpha}(G_1)=N,
\end{align}
where $N=\frac{\alpha n+s-1+\sqrt{(\alpha n+s-1)^{2}-4(\alpha^{2}sn-\alpha^{2}s^{2}+\alpha s^{2}-(1-\alpha)^{2}s\lfloor\frac{2}{3}s\rfloor-\alpha^{2}s+\alpha s-s)}}{2}$. We are to show $\lambda_{\alpha}(G_1)<n-2$. In terms of $n=\lfloor\frac{5}{3}s\rfloor+1$, we get
\begin{align}\label{eq:3.32}
(2(n-2)&-\alpha n-s+1)^{2}-(\alpha n+s-1)^{2}+4(\alpha^{2}sn-\alpha^{2}s^{2}+\alpha s^{2}-(1-\alpha)^{2}s\Big\lfloor\frac{2}{3}s\Big\rfloor-\alpha^{2}s+\alpha s-s)\nonumber\\
=&(4-4\alpha)n^{2}+(4\alpha^{2}s-4s+8\alpha-12)n-4\alpha^{2}s^{2}+4\alpha s^{2}-4(1-\alpha)^{2}s\Big\lfloor\frac{2}{3}s\Big\rfloor\nonumber\\
&-4\alpha^{2}s+4\alpha s+4s+8\nonumber\\
=&\left\{
\begin{array}{ll}
\frac{4}{9}((4-4\alpha)s^{2}+(9\alpha-15)s+9\alpha),&\mbox{if} \ s\equiv0 \ (\mbox{mod 3});\\
\frac{4}{9}((4-4\alpha)s^{2}+(17\alpha-23)s+5\alpha+10),&\mbox{if} \ s\equiv1 \ (\mbox{mod 3});\\
\frac{4}{9}((4-4\alpha)s^{2}+(13\alpha-19)s+8\alpha+4),&\mbox{if} \ s\equiv2 \ (\mbox{mod 3}).\\
\end{array}
\right.
\end{align}

\noindent{\bf Subcase 3.1.} $s\equiv0$ (mod 3).

We easily see $n=\frac{5}{3}s+1\geq25$, and so $s\geq15$. Write $\Phi_1(s)=(4-4\alpha)s^{2}+(9\alpha-15)s+9\alpha$. Since $0\leq\alpha<\frac{2}{3}$ and
$-\frac{9\alpha-15}{2(4-4\alpha)}<15\leq s$, we possess
\begin{align}\label{eq:3.33}
\Phi_1(s)\geq\Phi_1(15)=9(75-84\alpha)>0.
\end{align}

\noindent{\bf Subcase 3.2.} $s\equiv1$ (mod 3).

Obviously, $n=\frac{5s+1}{3}\geq25$, and so $s\geq16$. Let $\Phi_2(s)=(4-4\alpha)s^{2}+(17\alpha-23)s+5\alpha+10$. Since $0\leq\alpha<\frac{2}{3}$
and $-\frac{17\alpha-23}{2(4-4\alpha)}<16\leq s$, we infer
\begin{align}\label{eq:3.34}
\Phi_2(s)\geq\Phi_2(16)=9(74-83\alpha)>0.
\end{align}

\noindent{\bf Subcase 3.3.} $s\equiv2$ (mod 3).

Clearly, $n=\frac{5s+2}{3}\geq25$, and so $s\geq17$. Let $\Phi_3(s)=(4-4\alpha)s^{2}+(13\alpha-19)s+8\alpha+4$. Since $0\leq\alpha<\frac{2}{3}$ and
$-\frac{13\alpha-19}{2(4-4\alpha)}<17\leq s$, we deduce
\begin{align}\label{eq:3.35}
\Phi_3(s)\geq\Phi_3(17)=9(93-103\alpha)>0.
\end{align}

It follows from \eqref{eq:3.31}, \eqref{eq:3.32}, \eqref{eq:3.33}, \eqref{eq:3.34} and \eqref{eq:3.35} that $\lambda_{\alpha}(G_1)<n-2$. Together
with \eqref{eq:3.1} and \eqref{eq:3.4}, we conclude $\lambda_{\alpha}(G)\leq\lambda_{\alpha}(G_1)<n-2<\lambda_{\alpha}(K_1\vee(K_{n-2}\cup K_1))$,
which contradicts $\lambda_{\alpha}(G)\geq\lambda_{\alpha}(K_1\vee(K_{n-2}\cup K_1))$. This completes the proof of Theorem 1.1. \hfill $\Box$

\section*{Data availability statement}

My manuscript has no associated data.

\section*{Declaration of competing interest}

The authors declare that they have no conflicts of interest to this work.


\end{document}